\documentclass{article}
\usepackage{latexsym}
\usepackage[varumlaut]{yfonts}
\usepackage{epsfig}
\usepackage{amsmath}
\usepackage{amssymb}
\usepackage{amsthm}
\newtheorem{theorem}{Theorem}[section]
\newtheorem{corollary}[theorem]{Corollary}

\newtheorem{lemma}[theorem]{Lemma}

\theoremstyle{remark}

\theoremstyle{definition}
\newtheorem{definition}[theorem]{Definition}

\begin{document}

\title{Hessian potentials with parallel derivatives}

\author{Roland Hildebrand \thanks{%
LJK, Universit\'e Grenoble 1 / CNRS, 51 rue des Math\'ematiques, BP53, 38041 Grenoble cedex 09, France
({\tt roland.hildebrand@imag.fr}).}}

\maketitle

\begin{abstract}
Let $U \subset \mathbb A^n$ be an open subset of real affine space. We consider functions $F: U \to \mathbb R$ with non-degenerate Hessian such that the first or the third derivative of $F$ is parallel with respect to the Levi-Civita connection defined by the Hessian metric $F''$. In the former case the solutions are given precisely by the logarithmically homogeneous functions, while the latter case is closely linked to metrised Jordan algebras. Both conditions together are related to unital metrised Jordan algebras. Both conditions combined with convexity provide a local characterization of canonical barriers on symmetric cones.
\end{abstract}



\section{Introduction} \label{sec_intro}

Connections between different areas have since long been a major stimulant of progress in mathematics. One such connection is that between Jordan algebras and symmetric spaces. Convex symmetric cones, i.e., self-dual homogeneous convex cones, can be represented as the cones of squares of Euclidean Jordan algebras and are in one-to-one correspondence with these algebras \cite[Theorem 3]{Vinberg60}. This correspondence can be generalized to the non-convex case. A (not necessarily convex) symmetric cone is a conic domain carrying the structure of a symmetric space such that the transvections are linear maps. Non-degenerate such domains can be represented as the connection components of the identity of the set of invertible elements in a semi-simple Jordan algebra and are in one-to-one correspondence with these algebras \cite[Theorem 4]{Vinberg60}, see also \cite{Koecher99},\cite{Loos69a},\cite{Koecher70}. Hirzebruch \cite{Hirzebruch65} showed that the set of primitive idempotents of a simple Euclidean Jordan algebra is a compact connected Riemannian symmetric space of rank 1, and that any such space can be represented in this form. Geometric objects such as antipodal sets or the isometry group have an algebraic interpretation in terms of the Euclidean Jordan algebra. In \cite{Loos71} a one-to-one correspondence between compact Jordan triple systems and symmetric $R$-spaces was established. In \cite{Shima99} a one-to-one correspondence between semi-simple symmetric spaces with invariant projectively flat affine connection and central-simple Jordan algebras was established. If the symmetric space has a Riemannian metric, the corresponding Jordan algebra is Euclidean \cite{MizuharaShima99}.

Here we report another connection between Jordan algebras and geometry. We consider functions $F$ defined on a subset of an affine real space $\mathbb A^n$, with non-degenerate Hessian, such that the third derivative $F'''$ is parallel with respect to the Levi-Civita connection $\hat\nabla$ defined by the Hessian pseudo-metric $F''$. This condition can be written as a quasi-linear fourth order partial differential equation (PDE) on $F$. We show that the integrability condition of this PDE is equivalent to the condition that the commutative algebra defined by the difference tensor between the Levi-Civita connection $\hat\nabla$ and the canonical flat affine connection $\nabla$ on $\mathbb A^n$ is a Jordan algebra.

The construction applies also to Hessian manifolds. A {\it Hessian manifold} is a pseudo-Riemannian manifold equipped with a flat affine connection $\nabla$ such that the pseudo-metric $g$ can locally be expressed as the Hessian $\nabla^2F$ of a scalar function $F$, the {\it Hessian potential}. The Hessian potential is determined by $g$ up to an additive linear term. The parallelism condition $\hat\nabla(\nabla^3F) = 0$ considered in the previous paragraph is then equivalent to the condition $\hat\nabla\nabla g = 0$, which does not depend on the choice of the potential. We show (Theorem \ref{main3}) that the equivalence classes with respect to local isomorphism of Hessian manifolds satisfying this condition are in one-to-one correspondence with the isomorphism classes of metrised Jordan algebras, i.e., Jordan algebras equipped with a non-degenerate symmetric invariant bilinear form \cite{Bordemann97}. Any such manifold is a locally symmetric space.

We also consider the condition that the first derivative $F'$ is parallel with respect to the Levi-Civita connection defined by $F''$. We show (Theorem \ref{main1}) that a function $F$ satisfies this condition if and only if it is locally logarithmically homogeneous with respect to some central point in $\mathbb A^n$. This yields a differential-geometric description of logarithmic homogeneity.

We show (Theorem \ref{main13}) that given a Hessian manifold satisfying $\hat\nabla\nabla g = 0$, the existence of a Hessian potential satisfying $\hat\nabla\nabla F = 0$ is equivalent to unitality of the metrised Jordan algebra associated to the manifold. The conditions $\hat\nabla\nabla^3F = 0$, $\hat\nabla\nabla F = 0$ together with the convexity of $F$ characterize the class of canonical barriers on convex symmetric cones (Theorem \ref{main13convex}), which is extensively used in convex programming \cite{BenTalNemirovski01}.

Closely related to the present paper is the description of affine hypersurface immersions with parallel cubic form, which has been elaborated in \cite{Hildebrand12a},\cite{Hildebrand13a}.

\section{Notations and preliminaries}

Let $U \subset \mathbb A^n$ be a connected domain in the $n$-dimensional real affine space. We consider sufficiently smooth functions $F: U \to \mathbb R$ with non-degenerate Hessian. Such a function $F$ turns $U$ into a Hessian pseudo-Riemannian manifold. Let $\nabla$ be the canonical flat affine connection on $U$, $\hat\nabla$ the Levi-Civita connection of the Hessian pseudo-metric, and $K = \nabla - \hat\nabla$ the difference tensor. The difference tensor is a tensor of type $(1,2)$ which is symmetric in the two lower indices.

Denote the derivatives of $F$ with respect to the affine connection $\nabla$ by indices after a comma. Thus we have $\nabla_{\alpha}F = F_{,\alpha}$, $\nabla_{\alpha}\nabla_{\beta}F = F_{,\alpha\beta}$ etc. Denote the elements of the inverse of the Hessian $F'' = \nabla^2F$ by $F^{,\alpha\beta}$. We adopt the Einstein summation convention over repeating indices. Lowering and raising indices will be accomplished by contraction with the metric tensor $F_{,\alpha\beta}$ and its inverse $F^{,\alpha\beta}$, respectively.

In an affine coordinate system, the Christoffel symbols of the Levi-Civita connection $\hat\nabla$ have the form $\Gamma^{\gamma}_{\alpha\beta} = \frac12F_{,\alpha\beta\delta}F^{,\gamma\delta}$. Hence the difference tensor can be expressed by the derivatives of $F$ by
\begin{equation} \label{difference_tensor}
K^{\gamma}_{\alpha\beta} = -\Gamma^{\gamma}_{\alpha\beta} = -\frac12F_{,\alpha\beta\delta}F^{,\gamma\delta}.
\end{equation}

By $I$ we shall denote the identity matrix.

By definition of the Levi-Civita connection, the second derivative $F'' = \nabla^2F$ of $F$ is parallel with respect to $\hat\nabla$. In the next sections we shall consider the conditions that $F''' = \nabla^3F$ and $F' = \nabla F$ are parallel with respect to $\hat\nabla$.

\section{Parallel third derivative}

The covariant derivative of $F'''$ with respect to $\hat\nabla$ is given by
\begin{equation} \label{parallel_F3}
\hat\nabla_{\delta}F_{,\alpha\beta\gamma} = F_{,\alpha\beta\gamma\delta} - \frac12F^{,\rho\sigma}(F_{,\alpha\beta\rho}F_{,\gamma\sigma\delta} + F_{,\alpha\gamma\rho}F_{,\beta\sigma\delta} + F_{,\beta\gamma\rho}F_{,\alpha\sigma\delta}).
\end{equation}
Hence $F'''$ is parallel with respect to $\hat\nabla$ if and only if $F$ is a solution of the quasi-linear fourth order PDE
\begin{equation} \label{quasi_lin_PDE}
F_{,\alpha\beta\gamma\delta} = \frac12F^{,\rho\sigma}(F_{,\alpha\beta\rho}F_{,\gamma\delta\sigma} + F_{,\alpha\gamma\rho}F_{,\beta\delta\sigma} + F_{,\alpha\delta\rho}F_{,\beta\gamma\sigma}).
\end{equation}
Note that $F$ is a solution of \eqref{quasi_lin_PDE} if and only if $F+l$ is a solution, where $l$ is an arbitrary affine-linear function, i.e., a function satisfying $\nabla^2l = 0$. The functions $F$ and $F+l$ define also the same pseudo-metric $F''$ and the same difference tensor $K$.

Let us deduce the integrability condition of PDE \eqref{quasi_lin_PDE}. Introduce affine coordinates $x^{\alpha}$ on $U$. Assume that $F$ is of class $C^5$. Differentiating \eqref{quasi_lin_PDE} with respect to the coordinate $x^{\eta}$ and substituting the appearing fourth order derivatives of $F$ by the right-hand side of \eqref{quasi_lin_PDE}, we obtain after simplification
\begin{eqnarray*}
F_{,\alpha\beta\gamma\delta\eta} &=& \frac14F^{,\rho\sigma}F^{,\mu\nu}\left(F_{,\beta\eta\nu}F_{,\alpha\rho\mu}F_{,\gamma\delta\sigma} + F_{,\alpha\eta\mu}F_{,\rho\beta\nu}F_{,\gamma\delta\sigma} + F_{,\gamma\eta\nu}F_{,\alpha\rho\mu}F_{,\beta\delta\sigma} + F_{,\alpha\eta\mu}F_{,\rho\gamma\nu}F_{,\beta\delta\sigma} \right. \\
&& + F_{,\beta\eta\nu}F_{,\gamma\rho\mu}F_{,\alpha\delta\sigma} + F_{,\gamma\eta\mu}F_{,\rho\beta\nu}F_{,\alpha\delta\sigma} + F_{,\beta\eta\nu}F_{,\delta\rho\mu}F_{,\alpha\gamma\sigma} + F_{,\delta\eta\mu}F_{,\rho\beta\nu}F_{,\alpha\gamma\sigma} \\
&& \left. + F_{,\delta\eta\nu}F_{,\alpha\rho\mu}F_{,\beta\gamma\sigma} + F_{,\alpha\eta\mu}F_{,\rho\delta\nu}F_{,\beta\gamma\sigma} + F_{,\delta\eta\nu}F_{,\gamma\rho\mu}F_{,\alpha\beta\sigma} + F_{,\gamma\eta\mu}F_{,\rho\delta\nu}F_{,\alpha\beta\sigma}\right).
\end{eqnarray*}
The right-hand side must be symmetric in all 5 indices. Commuting the indices $\delta,\eta$ and equating the resulting expression with the original one we obtain
\begin{eqnarray*}
\lefteqn{F^{,\rho\sigma}F^{,\mu\nu}\left(F_{,\beta\eta\nu}F_{,\delta\rho\mu}F_{,\alpha\gamma\sigma} + F_{,\alpha\eta\mu}F_{,\rho\delta\nu}F_{,\beta\gamma\sigma} + F_{,\gamma\eta\mu}F_{,\rho\delta\nu}F_{,\alpha\beta\sigma} \right.} \\
&& \left. - F_{,\beta\delta\nu}F_{,\eta\rho\mu}F_{,\alpha\gamma\sigma} - F_{,\alpha\delta\mu}F_{,\rho\eta\nu}F_{,\beta\gamma\sigma} - F_{,\gamma\delta\mu}F_{,\rho\eta\nu}F_{,\alpha\beta\sigma}\right) = 0.
\end{eqnarray*}
Raising the index $\eta$, we get by virtue of \eqref{difference_tensor} the integrability condition
\[ K_{\alpha\mu}^{\eta}K_{\delta\rho}^{\mu}K_{\beta\gamma}^{\rho} + K_{\beta\mu}^{\eta}K_{\delta\rho}^{\mu}K_{\alpha\gamma}^{\rho} + K_{\gamma\mu}^{\eta}K_{\delta\rho}^{\mu}K_{\alpha\beta}^{\rho} = K_{\alpha\delta}^{\mu}K_{\rho\mu}^{\eta}K_{\beta\gamma}^{\rho} + K_{\beta\delta}^{\mu}K_{\rho\mu}^{\eta}K_{\alpha\gamma}^{\rho} + K_{\gamma\delta}^{\mu}K_{\rho\mu}^{\eta}K_{\alpha\beta}^{\rho}.
\]
This condition is satisfied if and only if $K_{\alpha\mu}^{\eta}K_{\delta\rho}^{\mu}K_{\beta\gamma}^{\rho}u^{\alpha}u^{\beta}u^{\gamma}v^{\delta} = K_{\alpha\delta}^{\mu}K_{\rho\mu}^{\eta}K_{\beta\gamma}^{\rho}u^{\alpha}u^{\beta}u^{\gamma}v^{\delta}$ for all tangent vectors fields $u,v$ on $U$, which can be written as
\begin{equation} \label{int_cond}
K(K(K(u,u),v),u) = K(K(u,v),K(u,u)).
\end{equation}

Consider an arbitrary point $y \in U$. The difference tensor $K$ defines a bilinear map $T_yU \times T_yU \to T_yU$ by $(u,v) \mapsto K(u,v)$. Equipped with this bilinear map, the tangent space $T_yU$ becomes an algebra $A_y$. We shall denote the multiplication in this algebra by $\bullet$, such that $u \bullet v = K(u,v)$. The left multiplication operator with the element $u$ will be denoted by $L_u$, such that $L_uv = u \bullet v$ for all $u,v$. Further, we define the positive powers of an element $u$ recursively by $u^1 = u$, $u^{k+1} = u \bullet u^k = L_u^ku$. If the algebra has a unit element $e$, then we put also $u^0 = e$. We need the following definitions.

{\definition Let $A$ be an algebra with multiplication $\bullet$. A bilinear form $\sigma$ on $A$ is called {\it invariant} if it satisfies the condition
\begin{equation} \label{3symmetry}
\sigma(u,v \bullet w) = \sigma(u \bullet v,w)
\end{equation}
for all $u,v,w \in A$. }

If $A$ is commutative and $\sigma$ is symmetric, then condition \eqref{3symmetry} is equivalent to the condition that the operator $L_v$ is self-adjoint with respect to $\sigma$ for all $v$.

{\definition \cite{Bordemann97} A {\it metrised algebra} is a pair $(A,\sigma)$ such that $A$ is an algebra and $\sigma$ is a non-degenerate symmetric invariant bilinear form on $A$. }

{\lemma \label{metrisedJ} Let $F: U \to \mathbb R$ be a $C^5$ solution of \eqref{quasi_lin_PDE} and $y \in U$ a point. Let $A_y$ be the algebra defined by the difference tensor $K$ on $T_yU$, and $\sigma_y$ the bilinear form defined on $T_yU$ by the Hessian pseudo-metric $F'' = \nabla^2F$. Then the pair $(A_y,\sigma_y)$ is a metrised Jordan algebra. }

\begin{proof}
The tensor $K$ is symmetric in the lower indices, and hence the multiplication $\bullet$ of the algebra $A_y$ is commutative. Condition \eqref{int_cond} becomes equivalent to the Jordan identity $u \bullet (u^2 \bullet v) = u^2 \bullet (u \bullet v)$, and $A_y$ is a Jordan algebra.

For arbitrary vectors $u,v,w \in A_y$ we have
\begin{eqnarray} \label{3sigma}
\sigma_y(u \bullet v,w) &=& F_{,\beta\gamma}K_{\delta\rho}^{\beta}u^{\delta}v^{\rho}w^{\gamma} = -\frac12F_{,\beta\gamma}F_{,\delta\rho\sigma}F^{,\sigma\beta}u^{\delta}v^{\rho}w^{\gamma} = -\frac12F_{,\delta\rho\gamma}u^{\delta}v^{\rho}w^{\gamma} \\ &=& -\frac12F_{,\beta\delta}u^{\delta}F_{,\rho\gamma\sigma}F^{,\sigma\beta}v^{\rho}w^{\gamma} =
F_{,\delta\beta}u^{\delta}K_{\rho\gamma}^{\beta}v^{\rho}w^{\gamma} = \sigma_y(u,v \bullet w). \nonumber
\end{eqnarray}
Here the second and fifth relation come from \eqref{difference_tensor}. Hence the form $\sigma_y$ satisfies \eqref{3symmetry} and is an invariant form. Finally, $\sigma_y$ is non-degenerate and symmetric because $F''$ is.
\end{proof}

Note that if $l$ is a linear function on $U$, then $F$ and $F+l$ define the same metrised Jordan algebra $(A_y,\sigma_y)$ on $T_yU$. Lemma \ref{metrisedJ} hence also applies to Hessian manifolds satisfying the condition $\hat\nabla\nabla g = 0$.

{\lemma \label{isomorphic} Let $F: U \to \mathbb R$ be a $C^5$ solution of \eqref{quasi_lin_PDE}, defined on a connected set $U$, and let $y,y' \in U$ be different points. Let $(A_y,\sigma_y)$,$(A_{y'},\sigma_{y'})$ be the metrised Jordan algebras defined on the tangent spaces $T_yU,T_{y'}U$ as in Lemma \ref{metrisedJ}. Then $(A_y,\sigma_y)$,$(A_{y'},\sigma_{y'})$ are isomorphic. }

\begin{proof}
Let $\gamma$ be a smooth path connecting the points $y,y'$. The parallel transport with respect to the Levi-Civita connection $\hat\nabla$ along $\gamma$ defines a non-degenerate linear map $J: T_yU \to T_{y'}U$. Now both $F''$ and $F'''$ are parallel with respect to $\hat\nabla$. Hence the difference tensor $K$ is also parallel. It then follows that $J$ is an isomorphism mapping $(A_y,\sigma_y)$ to $(A_{y'},\sigma_{y'})$.
\end{proof}

In particular, a closed path leading back to $y$ defines an automorphism of the metrised Jordan algebra $(A_y,\sigma_y)$.

\medskip

We have seen how a solution of \eqref{quasi_lin_PDE} defines a metrised Jordan algebra. We shall now consider the reverse direction.

{\lemma \label{algebraF} Let $(A,\sigma)$ be a metrised Jordan algebra. Then there exists a neighbourhood $U \subset A$ of zero such that the analytic function $F: U \to \mathbb R$ defined by
\begin{equation} \label{F_power_series}
F(x) = \sum_{k=2}^{\infty} \frac{(-1)^k}{k}\sigma(x,x^{k-1})
\end{equation}
is a solution of \eqref{quasi_lin_PDE}. }

\begin{proof}
First note that the expression $\sigma(x,x^{k-1})$ is a homogeneous polynomial of degree $k$ in the entries of $x$, and the right-hand side of \eqref{F_power_series} is an ordinary Taylor series. It is also easily seen that the convergence radius of the series is nonzero, and hence $F$ is defined on some neighbourhood $U \subset A$ of zero. On this neighbourhood $F$ is analytic. By possibly shrinking $U$, we shall also assume that the matrix $I+L_x$ is regular for all $x \in U$.

The partial derivative of $x^k$ in the direction $u$ is given by
\begin{eqnarray*} \label{grad_xk}
\nabla_u x^k &=& \nabla_u (L_x^{k-1}x) = \sum_{l=1}^{k-1} L_x^{l-1}L_uL_x^{k-1-l}x + L_x^{k-1}u = \sum_{l=1}^{k-1} L_x^{l-1}L_ux^{k-l} + L_x^{k-1}u \\ &=& \sum_{l=1}^k L_x^{l-1}L_{x^{k-l}}u,
\end{eqnarray*}
where $L_{x^0}$ is by convention the identity matrix. The derivative of $F$ is then given by
\begin{eqnarray} \label{derFu}
\nabla_u F &=& \sum_{k=2}^{\infty} \frac{(-1)^k}{k} \left( \sigma(\nabla_u x,x^{k-1}) + \sigma(x,\nabla_u x^{k-1}) \right) \nonumber\\ &=& \sum_{k=2}^{\infty} \frac{(-1)^k}{k} \left( \sigma(u,x^{k-1}) + \sum_{l=1}^{k-1} \sigma(x,L_x^{l-1}L_{x^{k-1-l}}u) \right) \nonumber\\ &=& \sum_{k=2}^{\infty} \frac{(-1)^k}{k} \left( \sigma(x^{k-1},u) + \sum_{l=1}^{k-1} \sigma(L_{x^{k-1-l}}L_x^{l-1}x,u) \right) = \sum_{k=2}^{\infty} (-1)^k \sigma(x^{k-1},u) \nonumber\\ &=& \sum_{k=1}^{\infty} (-1)^{k+1} \sigma(x^k,u) = \sigma((I+L_x)^{-1}x,u),
\end{eqnarray}
where the fourth equality comes from power-associativity of the Jordan algebra $A$ and all sums define analytic functions on $U$. Note that $I+L_x$ and its inverse are self-adjoint with respect to $\sigma$.

The next derivatives are given by
\begin{eqnarray*}
\nabla_v\nabla_uF &=& \sigma((I+L_x)^{-1}v,u) - \sigma((I+L_x)^{-1}L_v(I+L_x)^{-1}x,u), \\
\nabla_v^2\nabla_uF &=& -2\sigma((I+L_x)^{-1}L_v(I+L_x)^{-1}v,u) + 2\sigma((I+L_x)^{-1}L_v(I+L_x)^{-1}L_v(I+L_x)^{-1}x,u), \\
\nabla_v^3\nabla_uF &=& 6\sigma(((I+L_x)^{-1}L_v)^2(I+L_x)^{-1}v,u) - 6\sigma(((I+L_x)^{-1}L_v)^3(I+L_x)^{-1}x,u).
\end{eqnarray*}
At $x = 0$ we hence get
\begin{eqnarray*}
\nabla_v\nabla_uF &=& \sigma(v,u), \\
\nabla_v^2\nabla_uF &=& -2\sigma(v^2,u), \\
\nabla_v^4F &=& 6\sigma(v^3,v) = 6\sigma(v^2,v^2).
\end{eqnarray*}
It follows that \eqref{quasi_lin_PDE} is satisfied at $x = 0$.

For $x \in U$, we shall identify the tangent space $T_xU$ with $A$. For $w \in A$, define the vector field $X_w(x) = (I+L_x)w$ on $U$. Note that for arbitrary vector fields $V,Y$ on $U$ we have $\nabla_V\nabla_YX_w = \nabla_VL_Yw = L_{\nabla_VY}w = \nabla_{\nabla_VY}X_w$. Here, with a little abuse of notation, we denote by $L_Yw$ the vector field which at the point $x$ is given by $L_yw = y \bullet w$ with $y = Y(x)$. Since the vector field $X_w$ is affine-linear, the Lie derivative ${\cal L}_{X_w}$ commutes with the directional derivative $\nabla$.


We shall now compute the Lie derivative of \eqref{parallel_F3} with respect to the vector field $X_w$. By \eqref{derFu} we have
\[ {\cal L}_{X_w}F = \sigma((I+L_x)^{-1}x,(I+L_x)w) = \sigma(x,w),
\]
and the Lie derivative of $F$ is a linear function. It follows that ${\cal L}_{X_w}\nabla F = \nabla{\cal L}_{X_w}F$ is a constant 1-form, and ${\cal L}_{X_w}\nabla^kF = \nabla^{k-1}{\cal L}_{X_w}\nabla F = 0$ for every $k \geq 2$. Hence the Lie derivative of \eqref{parallel_F3} with respect to $X_w$ vanishes on $U$.

Since $I+L_x$ is regular, we have $\{ X_w(x) \,|\, w \in A \} = T_xU$ for all $x \in U$. If $U$ is connected, which we may assume without restriction of generality, then we thus have that \eqref{quasi_lin_PDE} is satisfied identically on $U$. This completes the proof.
\end{proof}

Lemmas \ref{metrisedJ} and \ref{algebraF} show how to construct a metrised Jordan algebra from a solution of \eqref{quasi_lin_PDE} and vice versa. We now show that the corresponding maps are the inverse of one another in the sense of the following lemmas.

{\lemma \label{two_directions} Let $F: U \to \mathbb R$ be a $C^5$ solution of \eqref{quasi_lin_PDE} and $y \in U$ be a point. Let $(A_y,\sigma_y)$ be the metrised Jordan algebra defined by $F$ as in Lemma \ref{metrisedJ}. Let $\tilde U \subset A_y = T_yU$ be the neighbourhood of zero and $\tilde F: \tilde U \to \mathbb R$ the solution of \eqref{quasi_lin_PDE} defined by $(A_y,\sigma_y)$ as in Lemma \ref{algebraF}.

Then there exists a neighbourhood $V \subset U \cap (y + \tilde U)$ of $y$ such that the difference $d(x) = F(x) - \tilde F(x - y)$ is affine-linear on $V$. }

\begin{proof}
The functions $F(x)$ and $\tilde F(x - y)$ are both defined on $V$ and are $C^5$ solutions of \eqref{quasi_lin_PDE}. We shall now compare the second and third derivatives of these functions at $x = y$. For vectors $u,v,w \in T_yV = A_y$ we have by definition of $\sigma_y$ and by \eqref{3sigma} that
\[ F''(u,v) = \sigma_y(u,v),\qquad F'''(u,v,w) = -2\sigma_y(u \bullet v,w).
\]
On the other hand, the quadratic and cubic terms in the Taylor series \eqref{F_power_series} yield
\[ \tilde F''(u,u) = \sigma_y(u,u),\qquad \tilde F'''(u,u,u) = -2\sigma_y(u,u \bullet u).
\]
Thus the second and third derivatives of $F(x)$ and $\tilde F(x - y)$ coincide at $x = y$.

Consider a ray $\gamma(t) = y + tz$ emanating from $y$. On this ray equation \eqref{quasi_lin_PDE} defines an ordinary differential equation (ODE) on the vector of second and third derivatives of $F$ and $\tilde F$, respectively. By the preceding paragraph, both ODEs have the same initial condition at $t = 0$. Since the second derivative is non-degenerate, the ODEs satisfy the conditions of the Picard-Lindel\"of theorem \cite{Lindeloef94} on the existence and uniqueness of the solution. Therefore the restriction of the second derivative $\nabla^2d$ to the ray $\gamma$ is identically zero on some interval $[0,T]$ with $T > 0$. The Lipschitz constant of the right-hand side of the ODE, which is involved in the proof of the Picard-Lindel\"of theorem and defines a strictly positive lower bound on $T$, is a continuous function of the direction $z$ of the ray. This bound then also depends continuously on $z$. By letting $z$ running through the unit sphere, it follows that there exists a neighbourhood of $y$ where $\nabla^2d$ identically vanishes. On this neighbourhood $d$ is an affine-linear function. This completes the proof.
\end{proof}

{\lemma \label{two_directions2} Let $(A,\sigma)$ be a metrised Jordan algebra, let $U \subset A$ be the neighbourhood of zero and $F: U \to \mathbb R$ the solution of \eqref{quasi_lin_PDE} defined by $(A,\sigma)$ as in Lemma \ref{algebraF}. Let $(A_0,\sigma_0)$ be the metrised Jordan algebra defined by $F$ at the point $y = 0$ as in Lemma \ref{metrisedJ}. Then, under identification of $A$ with $T_0U$, we have $(A_0,\sigma_0) = (A,\sigma)$. }

\begin{proof}
By \eqref{F_power_series} we have for arbitrary $u \in T_0U$ that
\[ F''(u,u) = \sigma(u,u),\qquad F'''(u,u,u) = -2\sigma(u,u \bullet u),
\]
where $\bullet$ denotes the multiplication in $A$. From the first relation it follows that $\sigma_0 = \sigma$. Since $A$ is commutative and $\sigma$ is a symmetric invariant form, it follows from the second relation that for all $u,v,w \in T_0U$ we have $F'''(u,v,w) = -2\sigma(u \bullet v,w)$. By \eqref{difference_tensor} we then get $K(u,v) = u \bullet v$, which proves also $A_0 = A$.
\end{proof}

From Lemma \ref{two_directions} we have also the following corollary.

{\corollary Let $F$ be a $C^4$ solution of \eqref{quasi_lin_PDE}. Then $F$ is analytic. }

\begin{proof}
If $F$ is $C^4$, then the right-hand side of \eqref{quasi_lin_PDE} is continuously differentiable. But then the left-hand side is continuously differentiable, and $F$ is actually $C^5$. By Lemma \ref{two_directions} $F$ then locally coincides with an analytic function. Hence $F$ is analytic.
\end{proof}

Summarizing the preceding results, we establish the following theorem.

{\theorem \label{main3} The equivalence classes with respect to local isomorphism of Hessian pseudo-metrics $g$ satisfying the condition $\hat\nabla\nabla g = 0$ are in one-to-one correspondence with the isomorphism classes of metrised Jordan algebras. }

\begin{proof}
Recall that a Hessian pseudo-metric $g$ is locally determined as the second derivative of a Hessian potential $F$ with respect to the flat affine connection $\nabla$. This potential is defined up to an additive affine-linear term. The Hessian pseudo-metric satisfies $\hat\nabla\nabla g = 0$ if and only if $F$ is a solution of \eqref{quasi_lin_PDE}.

By Lemma \ref{metrisedJ} every Hessian pseudo-metric $g$ satisfying the condition $\hat\nabla\nabla g = 0$ then determines metrised Jordan algebras on the tangent spaces of the Hessian manifold. By Lemma \ref{isomorphic} the isomorphism class of these metrised Jordan algebras is locally constant. By Lemmas \ref{algebraF} and \ref{two_directions} the metrised Jordan algebra in turn determines locally the Hessian pseudo-metric $g$. Hence non-isomorphic pseudo-metrics define non-isomorphic metrised Jordan algebras. Finally, by Lemma \ref{two_directions2} every metrised Jordan algebra can be produced in this way. This completes the proof.
\end{proof}

Note that a pseudo-Riemannian manifold is locally symmetric if and only if its Riemann curvature tensor is parallel with respect to the Levi-Civita connection, $\hat\nabla R = 0$. For Hessian pseudo-metrics $g = F''$, the curvature tensor is a quadratic function of the third derivative $F''' = \nabla g$ \cite[eq.~(1.7)]{Duistermaat01}. Thus the condition $\hat\nabla\nabla g = 0$ is a sufficient condition for local symmetry, and every solution of \eqref{quasi_lin_PDE} defines a locally symmetric Hessian pseudo-metric. It remains open whether every locally symmetric Hessian manifold satisfies the condition $\hat\nabla\nabla g = 0$.

\section{Parallel first derivative}

The covariant derivative of $F' = \nabla F$ with respect to $\hat\nabla$ is given by $\hat\nabla_{\beta}F_{,\alpha} = F_{,\alpha\beta} - \frac12F_{,\delta}F^{,\gamma\delta}F_{,\alpha\beta\gamma}$. Hence $F'$ is $\hat\nabla$-parallel if and only if
\begin{equation} \label{F1parallel}
F_{,\delta}F^{,\gamma\delta}F_{,\alpha\beta\gamma} = 2F_{,\alpha\beta}.
\end{equation}
Let $F: U \to \mathbb R$ be a solution of \eqref{F1parallel}. Define the vector field $e^{\gamma} = -F_{,\delta}F^{,\gamma\delta}$ on $U$. We then have
\begin{equation} \label{diff_e}
\nabla_{\alpha}e^{\gamma} = -F_{,\alpha\delta}F^{,\gamma\delta} + F_{,\delta}F^{,\gamma\rho}F_{,\rho\sigma\alpha}F^{,\sigma\delta} = -\delta^{\gamma}_{\alpha} + 2F^{,\gamma\rho}F_{,\rho\alpha} = \delta^{\gamma}_{\alpha},
\end{equation}
where $\delta^{\gamma}_{\alpha}$ is the Kronecker symbol.

Let $x^{\alpha}$ be an affine coordinate system on $\mathbb A^n$. By \eqref{diff_e} the vector field $e$ differs from the position vector field $x$ by a constant $c = x - e$. This difference distinguishes a point $c \in \mathbb A^n$, which we call the {\it center}. If $c \in U$, then $e$ vanishes at $c$, however, this condition is not necessary for the definition of the center. Let us shift the coordinate system in $\mathbb A^n$ such that $c = 0$, and the position vector field $x$ coincides with $e$. By definition of $e$ we then have
\begin{equation} \label{Fx2}
F_{,\delta} + F_{,\gamma\delta}x^{\gamma} = 0.
\end{equation}
Integrating, we obtain
\begin{equation} \label{Fx}
F_{,\gamma}x^{\gamma} = \nu,
\end{equation}
where $\nu \in \mathbb R$ is an integration constant. Integrating \eqref{Fx} along the rays emanating from $c$, we obtain
\begin{equation} \label{log_hom}
F(\alpha x) = \nu\log\alpha + F(x)
\end{equation}
for all $x \in U$ and $\alpha > 0$ such that the ray segment between $x$ and $\alpha x$ lies in $U$. This means that $F$ is locally logarithmically homogeneous with homogeneity parameter $\nu$.

On the other hand, let $F: U \to \mathbb R$ be a locally logarithmically homogeneous function with homogeneity parameter $\nu$ and with non-degenerate Hessian. Differentiating \eqref{log_hom} with respect to $\alpha$ at $\alpha = 1$ yields \eqref{Fx}. Differentiating \eqref{Fx} yields \eqref{Fx2}. Differentiating \eqref{Fx2} and eliminating $x$ by virtue of \eqref{Fx2} then gives back \eqref{F1parallel}.

We obtain the following result.

{\theorem \label{main1} Let $F: U \to \mathbb R$ be a $C^3$ function defined on some connected domain $U \subset \mathbb A^n$. Suppose that $F$ has a non-degenerate Hessian and denote by $\hat\nabla$ the Levi-Civita connection of the Hessian pseudo-metric $g = F''$. Then the first derivative $F'$ is $\hat\nabla$-parallel if and only if $F$ is locally logarithmically homogeneous with some homogeneity parameter $\nu$ with respect to some central point $c \in \mathbb A^n$. \qed }

\section{Parallel first and third derivatives}

In this section we consider the situation when both $\hat\nabla\nabla F = 0$ and $\hat\nabla\nabla^3 F = 0$ are satisfied.

{\lemma \label{i21} Let $F: U \to \mathbb R$ be a solution of \eqref{quasi_lin_PDE}, $y \in U$ a point, and $(A_y,\sigma_y)$ the metrised Jordan algebra defined by $F$ as in Lemma \ref{metrisedJ}. If $F$ in addition satisfies \eqref{F1parallel}, then the Jordan algebra $A_y$ possesses a unit element, which is given by $e^{\gamma} = -F_{,\delta}F^{,\gamma\delta}$. }

\begin{proof}
Raising the index $\beta$ in \eqref{F1parallel}, we obtain by \eqref{difference_tensor} that $-F_{,\delta}F^{,\gamma\delta}K^{\beta}_{\alpha\gamma} = \delta_{\alpha}^{\beta}$. The left-hand side of this equation defines the multiplication operator $L_e$ corresponding to the vector $e^{\gamma} = -F_{,\delta}F^{,\gamma\delta}$. The right-hand side is the identity operator on $T_yU$, and hence $e$ is a unit element of $A_y$.
\end{proof}

{\lemma \label{i12} Let $F: U \to \mathbb R$ be a solution of \eqref{quasi_lin_PDE}, $U$ connected, $y \in U$ a point, and $(A_y,\sigma_y)$ the metrised Jordan algebra defined by $F$ as in Lemma \ref{metrisedJ}. If $A_y$ possesses a unit element, then there exists an affine-linear function $l$ on $U$ such that $F-l$ is locally logarithmically homogeneous. }

\begin{proof}
If $A_y$ has a unit element, then for every $y' \in U$ the similarly defined Jordan algebra $A_{y'}$ has also a unit element, because by Lemma \ref{isomorphic} the algebras $A_y$ and $A_{y'}$ are isomorphic. Let $e$ be the vector field on $U$ defined by these unit elements. We then have
\begin{equation} \label{e_identity}
K^{\beta}_{\alpha\gamma}e^{\gamma} = \delta^{\beta}_{\alpha}.
\end{equation}
Note that the difference tensor $K$ as well as the Kronecker symbol are $\hat\nabla$-parallel. It follows that $K(u,\hat\nabla_ve) = 0$ for all vector fields $u,v$. Note that in a unital algebra $L_u = 0$ implies $u = 0$. Therefore we get $\hat\nabla e = 0$ identically on $U$. Writing this out, we get $\frac{\partial}{\partial x^{\delta}}e^{\gamma} + \Gamma^{\gamma}_{\beta\delta}e^{\beta} = \frac{\partial}{\partial x^{\delta}}e^{\gamma} - K^{\gamma}_{\beta\delta}e^{\beta} = 0$. But $e$ is the unit element, and hence $K^{\gamma}_{\beta\delta}e^{\beta} = \delta^{\gamma}_{\delta}$. As in the previous section, the vector field $e$ then differs from the position vector field $x$ by a constant, and we may choose the affine coordinate system in $\mathbb A^n$ such that $x = e$.

Lowering the index $\beta$ in \eqref{e_identity}, we get by virtue of \eqref{difference_tensor} that $-\frac12F_{,\alpha\beta\gamma}x^{\gamma} = F_{,\alpha\beta}$. Hence the tensor $\nabla^2(\nabla_xF)$ given by
\[ \frac{\partial^2F_{,\gamma}x^{\gamma}}{\partial x^{\alpha}\partial x^{\beta}} = F_{,\gamma\alpha\beta}x^{\gamma} + F_{,\gamma\alpha}\delta^{\gamma}_{\beta} + F_{,\gamma\beta}\delta^{\gamma}_{\alpha} = F_{,\gamma\alpha\beta}x^{\gamma} + 2F_{,\alpha\beta}
\]
vanishes on $U$. It follows that $\nabla_xF$ is an affine-linear function on $U$, which we will denote by $l$.

For every affine-linear function $\tilde l$ on $\mathbb A^n$ it holds that the difference $\nabla_x\tilde l - \tilde l$ is a constant. Therefore we have $\nabla_x(F-l) = \nu$ for some real number $\nu$. As in the previous section, it follows that $F-l$ is locally logarithmically homogeneous with homogeneity parameter $\nu$. This completes the proof.
\end{proof}

Combining the two lemmas, we obtain the following theorem.

{\theorem \label{main13} Let $F: U \to \mathbb R$ be a solution of \eqref{quasi_lin_PDE}, $U$ connected, $y \in U$ a point, and $(A_y,\sigma_y)$ the metrised Jordan algebra defined by $F$ as in Lemma \ref{metrisedJ}. Then the following are equivalent.

1. $A_y$ possesses a unit element.

2. There exists an affine-linear function $l$ on $U$ such that $F-l$ is locally logarithmically homogeneous. }

\begin{proof}
The implication 1 $\Rightarrow$ 2 is the assertion of Lemma \ref{i12}. The reverse implication follows from Lemma \ref{i21} by substituting $F-l$ for $F$ and noting that the metrised Jordan algebra $(A_y,\sigma_y)$ does not change.
\end{proof}

We shall now add the convexity of the Hessian potential $F$ as a condition.

{\lemma \label{13convex} Let $F: U \to \mathbb R$ be a convex solution of both \eqref{quasi_lin_PDE} and \eqref{F1parallel}, $y \in U$ a point, and $(A_y,\sigma_y)$ the metrised Jordan algebra defined by $F$ as in Lemma \ref{metrisedJ}. Then $A_y$ is a Euclidean Jordan algebra. }

\begin{proof}
We have to show that if $x_1^2 + \dots + x_k^2 = 0$ for some $x_1,\dots,x_k \in A_y$, then $x_j = 0$ for all $j = 1,\dots,k$. By convexity of $F$ the bilinear form $\sigma_y$ is positive definite. By Lemma \ref{i21} the algebra $A_y$ is unital. Let $e$ be the unit element. We have (cf.~\cite[Theorem VI.12]{Koecher99})
\[ 0 = \sigma_y\left(e,\sum_{j=1}^k x_j^2\right) = \sum_{j=1}^k \sigma_y(e,x_j \bullet x_j) = \sum_{j=1}^k \sigma_y(x_j,x_j).
\]
But $\sigma_y(x_j,x_j) \geq 0$ for all $j$, and $\sigma_y(x_j,x_j) = 0$ if and only if $x_j = 0$. This proves our claim.
\end{proof}

In a Euclidean Jordan algebra $A$, there exists for every element $x \in A$ a complete system of mutually orthogonal idempotents $e^1,\dots,e^m$ and distinct reals $\lambda_1,\dots,\lambda_m$ such that $x = \sum_{j=1}^m \lambda_j e^j$ \cite[Theorem 6]{JvNW34}. The numbers $\lambda_j$ are called the {\it eigenvalues} of $x$, and $d_j = tr L_{e^j}$ is their multiplicity. Clearly we have $x^k = \sum_{j=1}^m \lambda_j^k e^j$. The {\it determinant} of $x$ is defined as the product $\prod_{j=1}^m \lambda_j^{d_j}$. Let us compute the derivatives of the function $\log\det x$ at $x = e$. Fix a vector $u = \sum_{j=1}^{m'} \mu_j \tilde e^j$ and set $\tilde d_j = tr L_{\tilde e^j}$. Then we have $\log\det(e+tu) = \sum_{j=1}^{m'} \tilde d_j\log(1 + t\mu_j)$, and hence
\begin{equation} \label{der_det}
\nabla_u\log\det x|_{x = e} = \sum_{j=1}^{m'} \tilde d_j\mu_j = tr L_u,\qquad \nabla_u^2\log\det x|_{x = e} = -\sum_{j=1}^{m'} \tilde d_j\mu_j^2 = -tr L_{u^2}.
\end{equation}

{\lemma \label{simple} Let $A$ be a Euclidean Jordan algebra with unit element $e$, and let $\tau$ be the trace form on $A$, defined by $\tau(u,v) = tr\,L_{u \bullet v}$, where $\bullet$ is the multiplication in $A$. Then for every $x \in A$ with all eigenvalues in the open interval $(-1,1)$, we have
\[ -\log\det(e+x) = -tr L_x + \sum_{k=2}^{\infty} \frac{(-1)^k}{k}\tau(x,x^{k-1}).
\]
}

\begin{proof}
Let $x = \sum_{j=1}^m \lambda_j e^j$ be the eigenvalue decomposition of $x$ and $d_j$ the multiplicity of $\lambda_j$. We then have $e+x = \sum_{j=1}^m (1+\lambda_j) e^j$, and hence \begin{eqnarray*}
-\log\det(e+x) &=& -\sum_{j=1}^m d_j\log(1+\lambda_j) = \sum_{j=1}^m d_j\sum_{l=1}^{\infty}\frac{(-1)^l}{l}\lambda_j^l = \sum_{l=1}^{\infty}\frac{(-1)^l}{l}tr L_{x^l} \\ &=& -tr L_x + \sum_{k=2}^{\infty} \frac{(-1)^k}{k}\tau(x,x^{k-1}). \qedhere
\end{eqnarray*}
\end{proof}

Assume the conditions of Lemma \ref{13convex}, and let $\bullet$ denote the multiplication in $A_y$. By \cite[Theorems III.10,VI.12]{Koecher99} there exists a central element $z \in A_y$ (i.e., an element satisfying $z \bullet (u \bullet v) = u \bullet (z \bullet v)$ for all $u,v$) such that $\sigma_y(u,v) = tr\,L_{(z \bullet u) \bullet v}$ for all $u,v \in A_y$. By \cite[Corollary VI.5]{Koecher99} the algebra $A_y$ is semi-simple and can hence be decomposed as a direct sum $A_y = A^1_y \oplus \dots \oplus A^m_y$ of simple ideals. Let $e_1,\dots,e_m$ be the unit elements of these simple ideals. Then every central element $z \in A_y$ has a unique representation as a sum $z = \sum_{j=1}^m \alpha_je_j$, where $\alpha_j$ are real numbers \cite[p.46]{JvNW34}.

For a vector $u \in A_y$, let $u = \sum_{j=1}^m u_j$ be the decomposition of $u$ such that $u_j \in A^j_y$ for all $j$. Then we get $\sigma_y(u,v) = \sum_{j=1}^m \alpha_j tr\,L_{(e_j \bullet u) \bullet v} = \sum_{j=1}^m \alpha_j tr\,L_{u_j \bullet v_j}$ by the mutual orthogonality of the factors $A^j_y$. Note that the bilinear trace forms $\tau^j(u_j,v_j) = tr\,L_{u_j \bullet v_j}$ on $A^j_y$ are positive definite \cite[Theorem VI.12]{Koecher99}. In order for $\sigma_y = \sum_{j=1}^m \alpha_j \tau^j$ to be positive definite, we hence must have $\alpha_j > 0$ for all $j$. We get the following structural result.

{\lemma \label{structure13} Assume the conditions of Lemma \ref{13convex} and the notations of the previous paragraphs. Then in a neighbourhood of $y$ the function $F$ is given by
\begin{equation} \label{barrier_struct}
F(x) = -\sum_{j=1}^m \alpha_j \log\det\delta_j + F(y),
\end{equation}
where $\alpha_j > 0$ are some reals, $\delta = x - c$ is the difference between $x$ and some base point $c \in \mathbb A^n$, considered as an element of $A_y$, and $\det\delta_j$ denotes the determinant of $\delta_j$ in the simple factor $A^j_y$. }

\begin{proof}
Denote the difference $x - y$ by $\tilde\delta$ and consider it as a vector in $A_y$. By Lemma \ref{two_directions} there exists an affine-linear function $l$ on $\mathbb A^n$ such that $F(x) = l(x) + \sum_{k=2}^{\infty} \frac{(-1)^k}{k}\sigma_y(\tilde\delta,\tilde\delta^{k-1})$ in a neighbourhood of $y$. As we have shown above, there exist real numbers $\alpha_j > 0$ such that $\sigma_y(\tilde\delta,\tilde\delta^{k-1}) = \sum_{j=1}^m \alpha_j \tau^j(\tilde\delta_j,\tilde\delta_j^{k-1})$. From Lemma \ref{simple} it then follows that
\begin{equation} \label{F_intermed}
F(x) = l(x) + \sum_{j=1}^m \alpha_j \left( -\log\det(e_j+\tilde\delta_j) + tr L_{\tilde\delta_j} \right) = \tilde l(x) - \sum_{j=1}^m \alpha_j \log\det(e_j+\tilde\delta_j)
\end{equation}
with $\tilde l(x) = l(x) + tr L_{\tilde\delta}$ being another affine-linear function. By Lemma \ref{i21} the unit element of $A_y$ is given by $e^{\gamma} = -F_{,\delta}F^{,\gamma\delta}$, where the derivatives of $F$ are calculated at $y$. In other words, for every $u \in A_y$ we have $\sigma_y(e,u) = -\nabla_uF$. By \eqref{der_det} we have $\nabla_u F = \nabla_u \tilde l - \sum_{j=1}^m \alpha_j tr L_{u_j}$ at $x = y$. We then get
\[ -\nabla_u \tilde l + \sum_{j=1}^m \alpha_j tr L_{u_j} = \sigma_y(e,u) = \sum_{j=1}^m \alpha_j \tau^j(e_j,u_j) = \sum_{j=1}^m \alpha_j tr L_{u_j},
\]
and $\nabla_u \tilde l = 0$ for all $u$. Hence $\nabla \tilde l = 0$ at $y$. Since $\tilde l$ is affine-linear, it must be a constant. Setting $x = y$ in \eqref{F_intermed}, we get $\tilde l = F(y)$. Setting $c = y - e$, we get $\delta = e + \tilde\delta$. Hence \eqref{F_intermed} yields \eqref{barrier_struct}, which completes the proof.
\end{proof}

The class of functions defined by expression \eqref{barrier_struct} is well-known in convex conic optimization. These are the {\it canonical barriers} on convex symmetric cones \cite{BenTalNemirovski01}. Here the Euclidean Jordan algebra giving rise to the determinant in \eqref{barrier_struct} defines the cone, the base point $c$ in the argument $\delta = x - c$ locates the vertex of the cone, the weights $\alpha_j$ determine the contributions to the barrier of the determinants of the individual simple factors of the Jordan algebra, and the additive constant $F(y)$ is an offset which is irrelevant for the theory. Lemma \ref{structure13} has then the following consequence.

{\corollary \label{bar_struct} Let $F$ be a potential of a Hessian-Riemannian metric which satisfies both \eqref{quasi_lin_PDE} and \eqref{F1parallel}. Then $F$ is locally affinely isomorphic to a canonical barrier on a convex symmetric cone. \qed }

The converse is also true. Namely, a canonical barrier $F$ on a convex symmetric cone is by construction convex and logarithmically homogeneous. By Theorem \ref{main1} it then satisfies \eqref{F1parallel}. By Lemmas \ref{simple} and \ref{algebraF} the function $-\log\det\delta_j$ defined on any irreducible factor of the symmetric cone is a solution of \eqref{quasi_lin_PDE}. But then $F$ is also a solution of \eqref{quasi_lin_PDE}, because all derivatives of $F$ block-diagonalize in a coordinate system which is adapted to the decomposition of the symmetric cone in the product of its irreducible factors. We obtain the following result.

{\theorem \label{main13convex} The conditions \eqref{quasi_lin_PDE},\eqref{F1parallel}, and convexity locally characterize the class of canonical barriers on convex symmetric cones. \qed }

\bibliography{geometry,jordan,affine_geometry,algebra,misc,convexity}
\bibliographystyle{plain}

\end{document}